\documentclass[11pt]{article}
\usepackage[cp1251]{inputenc}
\usepackage[T2A]{fontenc}
\usepackage[english]{babel}
\usepackage{amsfonts}
\usepackage{amsmath}
\usepackage{amssymb}
\usepackage{graphicx}
\usepackage{multirow}
\usepackage[epsf]{}

\newtheorem{theorem}{Theorem}

\newtheorem{lemma}{Lemma}

\newcommand{\e}{\end{equation}}
\renewcommand{\b}{\begin{equation}}
\newcommand{\eps}{\varepsilon}
\newcommand{\proof}{{\bf Proof.\;}}

\title{Geometry of slow-fast Hamiltonian systems \\and Painlev\'e equations}
\author{L.M. Lerman, E.I. Yakovlev\\
\normalsize
Lobachevsky State University of Nizhny Novgorod, Russia}
\date{}
\begin{document}
\maketitle

\begin{abstract}
In the first part of the paper we introduce some geometric tools needed to
describe slow-fast Hamiltonian systems on smooth manifolds. We start with
a smooth Poisson bundle $p: M\to B$ of a regular (i.e. of constant rank) Poisson manifold $(M,\omega)$ over a
smooth symplectic manifold $(B,\lambda)$, the foliation into leaves of the bundle coincides
with the symplectic foliation generated by the Poisson structure on $M$. This defines a
singular symplectic structure $\Omega_{\varepsilon}=$ $\omega + \varepsilon^{-1}p^*\lambda$ on $M$
for any positive small $\varepsilon$, where $p^*\lambda$ is a lift of 2-form $\lambda$ on $M$.
Given a smooth Hamiltonian $H$ on $M$ one gets a slow-fast Hamiltonian system w.r.t.
$\Omega_{\varepsilon}$. We define a slow manifold $SM$ of
this system. Assuming $SM$ to be a smooth submanifold, we define a slow Hamiltonian
flow on $SM$. The second part of the paper deals with singularities of the
restriction of $p$ on $SM$ and their relations with the description of the
system near them. It appears, if $\dim M = 4,$ $\dim B = 2$ and Hamilton
function $H$ is generic, then behavior of the system near singularities of
the fold type is described in the principal approximation by the equation
Painlev\'e-I, but if a singular point is a cusp, then the related equation
is Painlev\'e-II. This fact for particular types of Hamiltonian systems
with one and a half degrees of freedom was discovered earlier by R.Haberman.
\end{abstract}

\section{Introduction}

Slow-fast Hamiltonian systems are ubiquitous in the applications in
different fields of science.
These  applications range from astrophysics, plasma physics and ocean hydrodynamics till molecular
dynamics. Usually these problems are given in the coordinate form,
moreover, in the form where a symplectic structure in the phase space is
standard (in Darboux coordinates). But there are cases where this form is
either nonstandard or the system under study is of a kind when its symplectic form has to
be found, in particular, when we deal with the system on a manifold.

It is our aim in this paper to present basic geometric tools to describe slow-fast Hamiltonian systems on manifolds,
that is in a coordinate-free way. For a general case this was done by V.I. Arnold \cite{Arn}. Recall that
a customary slow-fast dynamical system is defined by a system of differential equations
\begin{equation}\label{sf}
\eps\dot x = f(x,y,\eps),\;\dot y = g(x,y,\eps),\;(x,y)\in \mathbb R^m\times \mathbb R^n,
\end{equation}
depending on a small positive parameter $\eps$ (its positivity is needed to fix the positive direction of varying time $t$).
It is evident that $x$-variables in the region of the phase space where $f\ne 0$ change with the speed $\sim 1/\eps$ that is fast.
In comparison with them the change of $y$-variables is slow. Therefore variables $x$ are called to be fast and those of $y$ are slow ones.

With such the system two limiting systems usually connect whose properties influence on the dynamics of the slow-fast system
for a small $\eps.$ One of the limiting system is called to be fast or layer system and is derived in the following way.
One is introduced
a so-called fast time $\tau = t/\eps$, after that the system w.r.t. differentiating in $\tau$ gains the parameter $\eps$
in the r.h.s. of the second equation and looses it in the first equation, that is, the right hand sides become dependent in $\eps$
in a regular way
\begin{equation}\label{fs}
\frac{d x}{d\tau} = f(x,y,\eps),\;\frac{d y}{d\tau} = \eps g(x,y,\eps),\;(x,y)\in \mathbb R^m\times \mathbb R^n.
\end{equation}
Setting then $\eps =0$ we get the system, where $y$ are constants $y=y_0$ and they can be considered as parameters for the equations in $x$.
Sometimes these equations are called as layer equations. Because this system depends on parameters, it may pass through many bifurcations as
parameters change and this can be useful to find some special motions in the full system as $\eps > 0$ is small.

The slow equations are derived as follows. Let us formally set $\eps = 0$ in the system (\ref{sf}) and solve the equations
$f=0$ with respect to $x$ (where it is possible). The most natural case, when this can be done, is if matrix $f_x$
be invertible at the solution point in some domain where solutions for equations $f=0$ exist. Denote the related branch of solutions as $x=h(y)$
and insert it into the second equation instead of $x$. Then one gets a differential system w.r.t.
$y$ variables
$$
\dot y = g(h(y),y,0),
$$
which is called to be the slow system. The idea behind this construction
is as follows: if fast motions are directed to the slow manifold, then in
a small enough neighborhood of this manifold the motion of the full system
happens near this manifold and it is described in the first approximation by the
slow system.

Now the primary problem for the slow-fast systems is formulated as follows. Suppose we know something about the dynamics of both (slow
and fast) systems, for instance, some structure in the phase space composed from pieces of fast and slow motions.
Can we say anything about the dynamics of the full system for a small positive $\eps$ near this structure? There is a vast literature
devoted to the study of these systems, see, for instance, some of the references in \cite{Gucken}.

This set-up can be generalized to the case of manifolds in a free-coordinate manner
\cite{Arn}. Consider a smooth bundle $M\to B$ with a leaf $F$ being a smooth manifold and
assume a vertical vector field $v$ is given on $M$. The latter means
that any vector $v(x)$ is tangent to the leaf $F_b$ for any $x\in M$ and $b = p(x)\in B.$
In other words, every leaf $F_b$ of
$v$ is an invariant submanifold for this vector field. Let
$v_\eps$ be a smooth unfolding of $v = v_0$. Consider the set of
zeroes for vector field $v$, that is, one fixes a leaf $F_b$ then
on this smooth manifold $v$ generates a vector field $v^b$ and we consider its
zeroes (equilibria for this vector field). If the linearization operator of $v^b$ (along the leaf)
at zero $x$, being a linear operator $Dv^b_x: T_xF_b \to T_xF_b$ in
invariant linear subspace $V_x = T_xF_b$ of $T_xM$, has not zero eigenvalues,
then the set of zeroes is smoothly continued in $b$ for $b$ close to $b =
p(x)$. It is a consequence of the implicit function theorem. For this case
one gets a local section $z: B \to M,$ $p\circ z(b) = b$, which gives a smooth
submanifold $Z$ of dimension $dim B.$ One can define a vector field on $Z$
in the following way. Let us represent vector $v_\eps(x)$ in the unique way
as $v_\eps(x)= v^1_\eps(x) \oplus v^2_\eps(x)$, a
sum of two vectors of which $v^1_\eps(x)$ belongs to $V_x$ and $v^2_\eps(x)$ is in
$T_xZ$. Then vector $v^2_\eps(x)$ is of order $\eps$ in its norm, since
$v_\eps$ smoothly depends on $\eps$, and it is zero vector as $\eps = 0$.
Due to Arnold \cite{Arn} the vector field on $Z$ given as $(d/d\eps)(v^2_\eps)$ at $\eps = 0$ is
called to be slow vector field, in coordinate form it gives just what was
written above.

It is worth mentioning that one can call as {\em slow manifold} all set in
$M$ being the zero set for all vertical fields (for $\eps =0$). Generically, this set is a
smooth submanifold in $M$ but it can be tangent to leaves $F_b$ at some of its
points. In this case it is also possible sometimes to define a vector
field on $Z$ that can be called a slow vector field, but it is more
complicated problem intimately related with degenerations of the
projection of $p$ at these points (ranks of $Dp$ at these points, etc.) \cite{Arn}.

\section{Hamiltonian slow-fast systems}

Now we turn to Hamiltonian vector fields. It is well known, in order to define in an invariant way a Hamiltonian vector field,
the phase manifold $M$ has to be symplectic: a smooth nondegenerate closed
2-form $\Omega$ has to be given on $M$.
A slow-fast Hamiltonian system with a Hamiltonian $H(x,y,u,v,\eps)$ in coordinates
 $(x,y,u,v)=(x_1,\ldots,x_n,y_1,\ldots,y_n,u_1,\ldots,u_m,v_1,\ldots,v_m)$ has the form
\begin{equation}\label{sfH}
\begin{array}{l}
\displaystyle{\eps\dot{x}_i = \frac{\partial H}{\partial y_i},\;
\eps\dot{y}_i = -\frac{\partial H}{\partial x_i},\;i= 1,\ldots, n,}\\\\
\displaystyle{\dot{u}_j = \frac{\partial H}{\partial v_j},\;
\dot{v}_j = -\frac{\partial H}{\partial u_j},\;j= 1,\ldots, m.}
\end{array}
\end{equation}
It is worth mentioning that a symplectic structure for this system is
given by 2-form $\eps dx\wedge dy + du\wedge dv$ which is regularly
depends on $\eps$ and it degenerates into the Poisson co-symplectic structure at
$\eps =0$. Instead, if we introduce fast time $t/\eps = \tau$, then for the
transformed system one has 2-form $dx\wedge dy + \eps^{-1}du\wedge dv$
which depends singularly in $\eps.$

Let us note that the fast system here is Hamiltonian one (evidently) but the same is true for the slow system. Indeed,
if $x=p(u,v),$ $y = q(u,v)$ represent solutions of the system $H_x =0,$ $H_y =0$ (the set of zeroes for the fast system),
then the differential system
$$
\dot{u} = H_v,\;\dot{v} = - H_u,\; H(p(u,v),q(u,v),u,v,0)
$$
where {\em after differentiation} one needs to plug $x=p(u,v),$ $y = q(u,v)$ into Hamiltonian, is Hamiltonian with
the Hamilton function $h(u,v)= H(p(u,v),q(u,v),u,v,0),$ that is verified by
the direct differentiation.

From this set-up it is clear that the fast system is defined on a Poisson manifold with a regular (of the same rank at any point)
Poisson structure. Let us consider a smooth bundle $M\to B$ with a $C^\infty$-smooth connected Poisson manifold $(M,\omega)$ of
a constant rank with a co-symplectic structure $\omega$ \cite{Wein,Vais}. We assume that a symplectic foliation of $M$ into symplectic
leaves is just foliation into bundle leaves. Locally such symplectic foliation exists due to the Weinstein splitting theorem \cite{Wein}.
The dimension of these leaves, due to regularity and connectivity of $M$, is the integer $2n$ being the same for all
$M$ and it is just the rank of the Poisson bracket.

Now assume in addition that $B$ is a smooth symplectic manifold with a symplectic 2-form
$\lambda$ and consider a singular symplectic structure on $M$ generated by
2-form $\Omega_\eps = \omega + \eps^{-1}p^*\lambda$ where $p^*\lambda$ is a
pullback of $\lambda$ on $M$.

\begin{lemma}
For any positive $\eps$ 2-form $\Omega_\eps$ is nondegenerate, that is, $(M,\Omega_\eps)$ is a
symplectic manifold.
\end{lemma}
The proof of this lemma will be given below.

Let $H$ be a smooth function on the total space $M$ of the bundle (the Hamiltonian) and $X_{H}$ be its
Hamiltonian vector field. As is known \cite{Wein,Vais}, any leaf of the foliation is a smooth invariant symplectic manifold of
this vector field. This means that given a smooth Hamiltonian on $M$ one gets a family of reduced Hamiltonian vector fields on symplectic leaves,
that is, a family of Hamiltonian systems depending on some number of parameters $b\in B$, this number is equal to $2m$ where $2(n+m)$ is the dimension
of $M$.

Hamiltonian vector fields on the symplectic leaves are usually called to be fast or layer
systems. Thus, the foliation of $M$ into symplectic leaves is defined only
by the Poisson structure on $M$ but families of Hamiltonian vector fields
on leaves are defined by the function $H$.
Let us set $\eps =0$ and fix some symplectic leaf $F_b\subset M$. Suppose the related Hamiltonian vector field $X_{H}|_{F_b}$
to have a singular point $p$, $dH(\xi)|_p =0$ for all $\xi \in T_{p}F_b.$ If this singular point has not zero eigenvalues for its
linearization along the leaf $F_b$ at $p$, then the singular point persists in parameters $b$ and all close to $F_b$ symplectic leaves
have singular points for their Hamiltonian vector fields, these singular points depend smoothly on parameters. This gives locally
near $p$ a smooth submanifold of the dimension $2m$ that intersect each leaf close to $F_b$ at a unique point.
Thus we get in $M$ some smooth local submanifold $SM$. We call it {\em slow manifold} of $X_H.$

More generally, we come to the following construction.
Let $M$ be a smooth manifold of an even dimension $2n+2m$ and  a Poisson 2-form $\omega$ of a constant rank $2n$, be given on $M$. At every point $u\in M$ tangent space $T_uM$ contains $2n$-dimensional subspace $V_u\subset T_uM$ on which form $\omega$
is nondegenerate. Denote $N_u=(T_uM)^{\perp}$ the skew-orthogonal supplement
of $T_uM$ w.r.t. $\omega$. Then one has $\dim N_u=2m$ and $T_uM=V_u\oplus N_u$, maps $u\to V_u$ and $u\to N_u$
are smooth distributions $V$ and $N$ on $M$.

Due to the Darboux-Weinstein theorem \cite{Wein, Vais} distributions $V$ and $N$ are integrable. Their integral manifolds form smooth foliations
$F^V$ and $F^N$ on $M$. Moreover, restrictions of $\omega$ on leaves of $F^V$ are symplectic forms but on leaves of foliation
$F^N$ form $\omega$ is identically zero. Let us remark that the foliation of
$M$ into leaves of $V$ can be very complicated \cite{MR}.

Suppose the following conditions hold
\begin{itemize}
\item
there exists $2m$-dimensional manifold $B$ and a submersion $p: M\to B$ whose leaves coincide with the leaves of foliation $F^V$;
\item
there is a symplectic 2-form $\lambda$ on $B$.
\end{itemize}

For any $\eps\in\mathbb{R}$ small enough let us set
\begin{equation}\label{construction omega}
\Omega^{\varepsilon}=\varepsilon\omega+p^*\lambda.
\end{equation}
This defines a 2-form $\Omega^{\eps}$ on $M$.
\begin{lemma}\label{omega nondegenerate}
For all $\eps \ne 0$ form $\Omega^{\varepsilon}$ is symplectic.
\end{lemma}
\proof
2-forms $\omega$ and $\lambda$ are closed. Hence one has
$$
d\Omega^{\eps}=d(\eps\omega)+d(p^*\lambda)=\eps d\omega+p^* d\lambda=0.
$$
Consider a Darboux-Weinstein chart $(U,\varphi)$ for $\omega$ on $M$. Then for any point $u\in U$ the determinant of matrix
$$
A=(\omega(\partial_i^{\varphi}(u),\partial_j^{\varphi}(u))),\quad i,j=1,\dots,2n,
$$
is equal to $1$ and for all $r,s=1,\dots,2m$ and $\alpha,\beta=1,\dots,2n+2m$ we get
$$\omega(\partial_{2n+r}^{\varphi}(u),\partial_{\beta}^{\varphi}(u))=
\omega(\partial_{\alpha}^{\varphi}(u),\partial_{2n+s}^{\varphi}(u))=0.$$
On the other hand, the restriction $dp|_{N_u}:N_u\to T_aN$ is isomorphism, here $a=p(u)$. Therefore vectors
$$X_{2n+1}=dp(\partial_{2n+1}^{\varphi}(u)),\dots,X_{2n+2m}=dp(\partial_{2n+2m}^{\varphi}(u))$$
compose a basis of $T_aN$. Thus, matrix
$$
R=(\lambda(X_{2n+r},X_{2n+s})),\quad r,s=1,\dots,2m,
$$
is non-degenerate.

The last point is to remark that w.r.t. a holonomous basis $\partial_{1}^{\varphi}(u),\dots,\partial_{2n+2m}^{\varphi}(u)$
of the space $T_uM$ the matrix of the form $\omega$ has the form
$$C=
\begin{pmatrix}
\varepsilon A & 0\\
0             & R
\end{pmatrix},
$$
thus $\det{C}=\varepsilon^k\det{A}\det{R}=\varepsilon^k\det{R}\ne 0$.$\blacksquare$

\section{Symplectic submanifolds}

Consider now a $2m$-dimensional smooth submanifold $S\subset M$.
\begin{lemma}\label{omega on S}
If $S$ possesses the properties
\begin{itemize}
\item
at each point $v\in S$ submanifold $S\subset M$ intersects transversely the leaf through $v$ of foliation $F^V$,
\item
$S$ is compact,
\end{itemize}
then there is $\varepsilon_0 > 0$ such that for all $\varepsilon\in (0,\varepsilon_0)$ the restriction of the form
$\Omega^{\varepsilon}$ on $S$ is non-degenerate and hence is a symplectic form.
\end{lemma}
\proof
Consider first an arbitrary point $u\in M$ and its projection $a=p(u)$.
Since $p:M\to B$ is a submersion, there is a neighborhood $U'\subset M$ of $u$, a neighborhood $W'\subset B$ of $a$, a $2n$-dimensional smooth manifold $Q$ and diffeomorphism
$\psi:W'\times Q\to U'$ such that $p\circ\psi(b,q)= b$ for any $b\in W'$ and $q\in Q$.

For point $a$ there is a Darboux chart $(W_a,\theta)$ of $B$ w.r.t. symplectic form $\lambda$.
Without loss of generality one may assume that $W_a\subset W'$. Let us set $U_u=\psi(W_a\times Q)$ and consider $q\in Q$ such that $u=\psi(a,q)$.
Since $N$, $Q$ are manifolds then for points $a$ and $q$ there are neighborhoods $W_a^0$ and $Q^0$ which closures are compact
and belong to $W_a$ and $Q$, respectively. Moreover, the set $U_u^0=\psi(W_a^0\times Q^0)$ is a neighborhood of point $u$ and its
closure is also compact and belongs to $U_u$.

On $W_a$ holonomous vector fields $Y_{r}=\partial_{r}^{\theta}$, $r=1,\dots,2m$ are given. At any point of $W_a$ matrix
$(\lambda(Y_{r},Y_{s}))$, $r,s=1,\dots,2m$, has the canonical form. Therefore one has
\begin{equation}\label{det of nu}
\det(\lambda(Y_{r},Y_{s}))\equiv 1.
\end{equation}
In accordance to the first condition of the lemma, for any point $v\in S\cap U_u$ one has $T_vM=V_v\oplus T_vS$.
Consequently, the restriction $dp|_{T_vS}:T_vS\to T_{p(v)}B$ is isomorphism. Setting $Y_r^*(v)=(dp|_{T_vS})^{-1}(Y_r(p(v))$
for all $v\in S\cap U_u$ and $r=1,\dots,2m$ we get smooth vector fields $Y_1^*,\dots,Y_{2m}^*$ on $S\cap U_u$.
At any point $v\in S\cap U_u$ vectors $Y_1^*(v),\dots,Y_{2m}^*(v)$ make up a basis of the tangent space $T_vS$.

Due to (\ref{construction omega}) we get
$$
\Omega^{\varepsilon}(Y_r^*,Y_s^*)=\varepsilon\omega(Y_r^*,Y_s^*)+\lambda(Y_r,Y_s)
$$
for all $r,s=1,\dots,2m$. Thus for matrix $D=(\Omega^{\varepsilon}(Y_r^*,Y_s^*))$ there is
\begin{equation}\label{det of omega on S}
\det{D}=f_{2m}\varepsilon^{2m}+\dots+f_1\varepsilon+f_0,
\end{equation}
where for any $t=0,1,\dots,2m$ coefficient $f_t$ is the sum of all determinants for which $t$ rows coincide with the related rows of matrix
$(\omega(Y_r^*,Y_s^*))$, but other rows belong to matrix $(\lambda(Y_r,Y_s))$.
It follows from here that $f_t:S\cap U_u\to\mathbb{R}$ are smooth functions and due to identity (\ref{det of nu}) one has
\begin{equation}\label{f_0}
f_0=\det(\lambda(Y_{r},Y_{s}))\equiv1.
\end{equation}

Because the closure of the set $S\cap U_u^0$ is compact and belongs to $S\cap U_u$, functions $f_t$ are bounded on $S\cap U_u^0$.
Hence, there is $\varepsilon_u>0$ such that
\begin{equation}\label{f_k}
|f_{2m}\varepsilon^{2m}+\dots+f_1\varepsilon|<1
\end{equation}
on $S\cap U_u^0$ for any $\varepsilon\in(0,\varepsilon_u)$.

A collection $\mathcal{U}=\{U_u^0|u\in S\}$ is a cover for manifold $S$ by open sets. Since $S$ is
compact we conclude $S\subset U_{u_1}^0\cup\dots\cup U_{u_l}^0$ for some finite set of points $u_1,\dots,u_l\in S$.
Let us set $\varepsilon_0=\min\{\varepsilon_{u_1},\dots,\varepsilon_{u_l}\}$. Then for any
$\varepsilon\in(0,\varepsilon_0)$ the inequality (\ref{f_k}) is valid on every set from $S\cap U_{u_1}^0,\dots,S\cap U_{u_l}^0$.

If now $v$ is any point of $S$ there exists an $i\in\{1,\dots,l\}$ such that $v\in S\cap U_{u_i}^0$.
Here for any $\varepsilon\in(0,\varepsilon_0)$ it follows from (\ref{det of omega on S}), (\ref{f_0}) and (\ref{f_k}) that $\det{D}\ne 0$.
This implies 2-form $\Omega^{\varepsilon}$ be non-degenerate on the tangent space  $T_vS$.$\blacksquare$

\section{Slow manifold and nearby orbit behavior}

Henceforth we assume a Poisson bundle $p: M \to B$ be given where $M$ be a
$C^\infty$-smooth Poisson manifold with 2-form $\omega$, $B$ be a smooth symplectic manifold
with its 2-form $\lambda$ and $p$ be a smooth bundle map whose leaves define
symplectic foliation of $M$ given by $\omega.$
Suppose a smooth function $H$ on $M$ be given. The set of all zeroes for all vertical (fast) Hamiltonian vector fields generated by $H$ forms
a subset in $M$ being generically a smooth submanifold $SM$ of dimension $2m =
\mbox{dim\;}B$. We assume this is the case and call $SM$ to be the {\em slow manifold} of the vector field
$X_H$. When restricted on $SM$ the related map $p_r: SM \to B$ may contain in $SM$ points of two types: regular and
singular. By a {\em regular} point $s\in SM$ one understands such that $rank Dp_r(s) = 2m = dim
B$. This implies $p_r$ be a diffeomorphism near $s$. To the contrast, a point $s\in
SM$ is {\em singular}, if rank of $Dp_r$ at $s$ is lesser than $2m$.

Another characterization of regular and singular points appeals to a type of related equilibria for the fast Hamiltonian system
on the symplectic leaf containing $s$. The point $s$ is regular if the fast vector field has a simple equilibrium for the related fast
Hamiltonian vector field, i.e. such an equilibrium on the related leaf of the symplectic foliation that has not zero
eigenvalues. To the contrast, for a singular point of $s\in SM$ the
equilibrium on the symplectic leaf through $s$ for the fast Hamiltonian
vector field is degenerate, that is it has a zero eigenvalue. All this
will be seen below in local coordinates though it is possible to show it
in a coordinate-free way.

\subsection{Near regular points of $SM$}

Manifold $SM$ near its regular point $s$ can be represented by the implicit function theorem as a graph of a
smooth section $z: U \to M$, $p(s)\in U\subset B,$ $p\circ z = id_U$.
Due to lemma \ref{omega on S}, such a compact piece of $SM$ is a symplectic
submanifold w.r.t. the restriction of 2-form $\Omega_{\eps}={\eps}^{-1}\Omega^{\eps}$ to $SM$.
Hence one can define a slow Hamiltonian vector field on it generated by function
$H$. In the same way one can consider the case when function $H$ itself depends smoothly on a parameter $\eps$
on $M$. Let $X_H$ be the Hamiltonian vector field on $M$ w.r.t. 2-form
$\Omega_\eps$ generated by $H$. Denote $H^S$ the restriction of $H$ to $SM$ and consider a
Hamiltonian vector field on $SM$ with the Hamiltonian $H^S$ w.r.t. the restriction of 2-form $\Omega_\eps$ on $SM$.
This vector field is of the order $\eps$, hence there is a limit for $\eps^{-1}H^S$ as $\eps \to 0.$ This
is what is called to be {\em slow Hamiltonian vector field} on $SM$.

It is an interesting problem to understand the orbit behavior of the full system (for small $\eps > 0$) within a small neighborhood
of a compact piece of regular points in $SM$. This question is very hard in the general
set-up. Nevertheless, there is a rather simple general case to examine, if one assumes
the hyperbolicity of this piece of $SM$. Suppose for a
piece of $SM$ fast systems have equilibria at $SM\cap F_b$ (for $\eps = 0$) without zero real parts
(hyperbolic equilibria in the common terminology (see, for instance, \cite{SSTC,Meiss}). Then this smooth
submanifold of the vector field $X_H$ on $M$ is
a normally hyperbolic invariant manifold and results of \cite{Fenichel,HPS} are applicable.
They say that for $\eps > 0$ small enough there is a true invariant manifold
in a $O(\eps)$-neighborhood of that piece of $SM$. For the full system this invariant manifold is
hyperbolic and possesses its stable and unstable local smooth invariant
manifolds. The restriction of $X_{H_\eps}$ on this slow manifold can
be any possible Hamiltonian system with $m$ degrees of freedom. If we add to
this picture the extended  structure of the fast system along with their
bifurcations w.r.t. slow variables (parameters of the fast system), then
one can say a lot on the full system as $\eps$ is positive small. This
behavior is the topics of the averaging theory, theory of adiabatic
invariants, etc., see, for instance, \cite{Verhulst,Neishtadt}.

Much more subtle problem is to understand the local dynamics of the full system near $SM$
when the fast dynamics possesses elliptic equilibria at points of $SM$, this corresponds to a one degree of freedom fast systems.
Nonetheless, one can present some details of this picture when we deal
with a real analytic case (manifolds and Hamiltonian). Then results of
\cite{GL1} can be applied. For this case $SM$ was called in \cite{GL1} to be an {\em almost invariant elliptic} slow
manifold. At $\eps = 0$ near a piece of the almost elliptic slow manifold  one can introduce a coordinate
frame where this slow manifold corresponds to a zero section of the bundle $M\to B$. Then the main result of \cite{GL1}
is applicable
for the case when the fast system is two dimensional but the slow system can be of any finite dimension
and Hamiltonian is analytic in a neighborhood of $SM$. The result says on the
exponentially exact (with respect small parameter $\eps$) reduction of the Hamiltonian to the form where fast
variables $(x,y)$ enter locally near $SM$ only in the combination $I=(x^2+y^2)/2$. Thus, up
to the exponentially small error the system has an additional integral
$I$. In particular, if the slow system is also two dimensional, this gives
up to exponentially small error an integrable system in a small neighborhood of $SM.$ All this helps
a lot when one is of interest in the dynamics within this neighborhood, see,
for instance, \cite{GL2}. The case when fast system has more degrees of
freedom and the related equilibria are multi-dimensional elliptic ones is
more hard and no results are known to the date. The case of fast
equilibria with eigenvalues in the complex plane lying both on the imaginary axis and
out of it is even lesser explored.

\subsection{Near a singular point of $SM$}

At a singular point $s\in SM$ submanifold $SM$ is tangent to the related leaf of the symplectic
foliation. More precisely, the rank of the restriction of $Dp$ to the tangent plane
$T_sSM$ is less than the maximal $2m$. This means $Ker Dp|_{SM} \ne \emptyset$.
If we choose some Darboux-Weinstein
coordinate chart $(x,y,u,v)$ near $s$, then the Poisson form in these coordinates is given as $\Omega_0 = dx\wedge dy$ and symplectic leaves
are given as $(u,v)=(u_0,v_0)$. Let $H(x,y,u,v)$ be a smooth Hamilton function written in these coordinates, $dH \ne 0$ at some point $s$.
The Hamiltonian vector field near $s$ is given as
$$
\dot x = H_y,\; \dot y = - H_x,\; \dot u =0, \;\dot v =0.
$$
The condition point $s$ be in $SM$ (a singular point for the fast vector field) means $H_y(s)=
H_x(s)=0$ and if this point is a singular point for projection $p$, then
this implies
$$
\det\begin{pmatrix}\frac{\partial^2 H}{\partial (x,y)^2}\end{pmatrix}|_s=0,
$$
otherwise the system $H_y = 0,\;H_x=0$ can be solved w.r.t. $x,y$ and points in $SM$ be regular.
Thus, the fast vector field at point $s$ has zero eigenvalue, thus its eigenspace is invariant w.r.t. the
linearization of the fast vector field at $s$. This is equivalent to
the condition that $Dp$ restricted on $T_sSM$ has nonzero kernel and hence its rank is lesser than $2m$.

The types of degenerations for the mappings of one smooth manifold to another one
(for our case it is $p_r:SM \to B$) are studied by the singularity theory of smooth manifolds
\cite{Whitney,Arn_GZ}. When the dimension of $B$ (and $M$) are large,
these degenerations can be very complicated. Keeping this in mind, we consider below only
the simplest case of one fast and one slow degrees of freedom. For this
case $SM$ and $B$ have dimension 2 and we get the mapping from one two
dimensional smooth manifold to another smooth two dimensional one. In such cases degenerations can be generically
only of two types: folds and cusps \cite{Whitney}. We shall show that near a fold point the
system can be reduced to the case of a family of slow varying Hamiltonian systems with the
Hamiltonians $H(x,y,\eps t,c)$ with scalar $x,y$ and positive small parameter
$\eps,$ $c$ is also a small parameter. For a fixed $c$ the slow manifold of this one and a half degree of freedom system
is a slow curve and
the fold point corresponds to the quadratic tangency of the slow curve with
the related two-dimensional symplectic leaf. The fast (frozen) Hamiltonian system with one degree of freedom has the equilibrium,
corresponding to  the tangency point, being generically a parabolic equilibrium point with a double non-semisimple zero
eigenvalue. The local orbit behavior for these fast system do not change when $c$ varies near a fixed $c_0$.

For the case of a cusp we get again a slow fast Hamiltonian system with two dimensional slow manifold but fast Hamiltonian
system has at $\eps =0$ on the related symplectic leaf an equilibrium of the type of
degenerate saddle or degenerate elliptic point (both are of codimension
2). In this case it is also possible to reduce the system to the case of a family of slow
varying nonautonomous Hamiltonian systems but the behavior of these
systems near $s$ depends essentially on parameter $c$.

Henceforth, we consider only the case of one slow and one fast degrees of freedom Hamiltonian systems,
that is $M$ be a smooth 4-dimensional regular Poisson manifold of rank 2 with 2-form $\omega$ and
$B$ be a smooth symplectic 2-dimensional manifold with a symplectic 2-form
$\lambda$, $p:M \to B$ be a Poisson bundle with leaves $F_b$ which coincide with symplectic fibration generated by
$\omega$. The symplectic structure on $M$ is given by 2-form $\Omega_\eps = \omega + \eps^{-1}p^*\lambda.$

\section{Folds for the slow manifold projection}

Let a smooth Hamilton function $H$ on $M$ be given (here and below it is assumed $C^\infty$-smoothness).
We suppose $H$ to be non-degenerate within a neighborhood
of a point $s$ we work: $dH(s) \ne 0$. Then levels $H=c$ are smooth 3-dimensional disks within this neighborhood.
Since the consideration is local, we can work in
Darboux-Weinstein coordinates near point $s$, hence we suppose 2-form $\Omega_\eps$
be written as $\Omega_\eps = dx\wedge dy + \eps^{-1}du\wedge dv$ with
fast variables $x,y$ and slow variables $u,v$. The related Poisson manifold is endowed locally with 2-form $\omega = dx\wedge
dy$, its symplectic leaves are given as $(u,v)=(u_0,v_0)\in U$, $U\subset B$ is a disk with coordinates $(u,v).$
The restriction of $H$ on a symplectic leaf is the function $H(x,y,u_0,v_0)$ and the orbit
foliation for this fast vector field on the leaf
is given in fact by level lines $H(x,y,u_0,v_0) = c.$
Without loss of generality we assume $s$ being the origin of the coordinate frame,
$s=(0,0,0,0).$ Henceforth we assume $s$ to be an equilibrium of the fast
Hamiltonian system on the leaf $(0,0):$ $H_y(0,0,0,0)= H_x(0,0,0,0)=0.$

The equilibrium $s$ for the fast Hamiltonian vector field
$$
\dot x = \frac{\partial H}{\partial y},\;\dot y = -\frac{\partial H}{\partial x},\;u,v\,\mbox{are\; parameters},
$$
varies in parameters $(u_0,v_0)$ and its unfolding in parameters gives a local piece of slow manifold
$SM$ near point $s=(0,0,0,0)$ in the space $(x,y,u,v).$ Locally near $s$ the slow manifold is indeed
a smooth 2-dimensional disk, if rank of the matrix
\begin{equation}\label{2X2}
\begin{pmatrix}H_{xy}&H_{yy}&H_{uy}&H_{vy}\\H_{xx}&H_{yx}&H_{ux}&H_{vx} \end{pmatrix}
\end{equation}
at $s$ is 2. We suppose this is the case. This is a genericity condition on the function $H$.
Then we may consider the restriction of the projection map $p$ on  $SM$, this generates
the mapping $p_r$ of two smooth 2-dimensional manifolds $p_r: SM \to B$. If the
inequality $\Delta = H_{xy}^2 - H_{yy}H_{xx} \ne 0$ holds at $s$ (i.e. $s$ is a regular point of $SM$
and simultaneously is a simple equilibrium of the related fast vector field, that is without zero eigenvalues),
then by the implicit function theorem the
set of solutions for the system $H_y =0,\;H_x=0$ near $s$ is expressed as $x = f(u,v),$ $y=g(u,v)$ and hence
locally it is a section of the bundle $p: M\to B$ and $Dp_r$ does not degenerate on this set in
some neighborhood of $s$: $p_r$ is a diffeomorpism. Also, we get regular points for $p_r$ near $s\in
SM$ in the sense of singularity theory \cite{Whitney,Arn_GZ}.

Thus, the degeneration only happens if $\Delta(s) = 0$. This equality is equivalent to a condition the fast Hamiltonian vector
field on the related symplectic leaf $F_b,$ $b=p(s),$
to have a degenerate equilibrium at point $s$: it possesses zero eigenvalue (by necessity being double) for the linearization at
$s$. Another characterization of such point is that it is a {\em singular} point
of mapping $p_r$: rank of this mapping at $s$ is lesser than 2. For the
further goals one needs to distinguish singular points of a general type
that are possible for smooth maps of one 2-dimensional manifold to another one.

According to \cite{Whitney}, a singular point $q$ of a $C^2$-smooth map $F: U\to V$ of two open domains in
smooth 2-dimensional manifolds is {\em
good}, if the function $J=\det\,DF$ vanishes at $q$ but its differential $dJ$ is nondegenerate at this point.
In a neighborhood of a good singular
point $q$ there is a smooth curve of other singular points continuing $q$ \cite{Whitney}. Let
$\varphi(\tau)$ be a smoothly parametrized curve of singular points through a good singular point $q$ and $\tau =0$ corresponds to $q$.
We shall use below the notation $A^\top$ for the transpose matrix of any matrix $A$.

According to \cite{Whitney}, a good singular point $q$ is called to be a {\em fold} point if
$dF(\varphi'(0))\ne (0,0)^\top$, and it is called to be a {\em cusp} point if at $q$ one has
$dF(\varphi'(0))= (0,0)^\top$ but $d^2(F\circ \varphi)/d\tau^2|_q \ne (0,0)^\top.$ It is worth remarking
that if $q$ is a fold, then for any nearby point on the curve $\varphi(\tau)$ there
is a unique (up to a constant) nonzero vector $\xi$ in the tangent space of the related point such that
$DF(\xi)=0.$ The direction spanned by $\xi$ is transverse to the tangent direction to the singular
curve at the fold point, since $DF(\varphi'(0))\ne (0,0)^\top.$ These transverse
directions form a smooth transverse direction field on the curve if this curve is at
least $C^3$-smooth.

Now let us return to the set of points in $SM$ near $s$ where  $p_r$ degenerates. To be
precise, we assume $H_{yy}H_{ux} - H_{xy}H_{uy} \ne 0$ at $s$, then $SM$ near
$s$ is represented as a graph of $y=f(x,v),$ $u=g(x,v)$, $f(0,0)=g(0,0)=0$ (we can always assume this is the case, otherwise, one can
achieve this by re-ordering slow or/and fast variables). The unique case when rank equals 2 but the only nonzero minor is
$H_{yu}H_{vx} - H_{xu}H_{yv}$, but all other five are zeroth, would indicate on a too degenerate case and we do not consider it
below (recall that we have assumed rank be 2 for matrix (\ref{2X2}). Indeed, in that case the following lemma is valid
\begin{lemma}\label{matrix 2x4}
Suppose matrix
$$A=
\begin{pmatrix}
a_{11} & a_{12} & a_{13} & a_{14}\\
a_{21} & a_{22} & a_{23} & a_{24}
\end{pmatrix}
$$
possesses the properties:
\begin{itemize}
\item
$$\det
\begin{pmatrix}
a_{13} & a_{14}\\
a_{23} & a_{24}
\end{pmatrix}\ne 0,
$$
\item all other minors of the second order for matrix $A$ vanish.
\end{itemize}
Then the following equalities hold $a_{11} = a_{12} = a_{21} = a_{22} = 0$.
\end{lemma}
\proof
Suppose the assertion of Lemma is false. Then up to re-enumeration of rows and first two columns one may regard  $a_{11}\ne 0$.
By condition, all minors of the second order for matrix
$$
\begin{pmatrix}
a_{11} & a_{12} & a_{13}\\
a_{21} & a_{22} & a_{23}
\end{pmatrix}
$$
vanish. This means that its rows are linear dependent. But by assumption the first row is nonzero vector of  $\mathbb{R}^3$.
Hence, there is $\kappa_1\in\mathbb{R}$ such that $(a_{21}, a_{22}, a_{23})=\kappa_1(a_{11}, a_{12}, a_{13})$.

Analogously, if all second order minors of matrix
$$
\begin{pmatrix}
a_{11} & a_{12} & a_{14}\\
a_{21} & a_{22} & a_{24}
\end{pmatrix}
$$
vanish, then from $a_{11}\ne 0$ follows the existence of a number $\kappa_2\in \mathbb{R}$ which satisfies the equality
$(a_{21}, a_{22}, a_{24})=\kappa_2(a_{11}, a_{12}, a_{14})$. But then one has $\kappa_1a_{11}=a_{21}=\kappa_2a_{11}$, from here
equalities $a_{11}(\kappa_1-\kappa_2)=0$ and $\kappa_1=\kappa_2$ follow.
Thus, rows of matrix $A$ are linear dependent that contradicts to its first property of the Lemma
conditions.$\blacksquare$

For the case under consideration the coordinate representation for mapping $p_r$ is the following $p_r: (x,v)\to
(u=g(x,v),v)$. Jacobian of this mapping is the matrix
$$
P=Dp_r =\begin{pmatrix}g_x&g_v\\0&1\end{pmatrix},
$$
whose rank at $(0,0)$ is 2, if $g_x(0,0)\ne 0$ (the point is regular), and it is 1, if $g_x(0,0)= 0$ (the point is singular).
The singular point $(0,0)$ is good, if $g_x(0,0)=0$ and $g_{xx}(0,0)\ne 0$ or $g_{xv}(0,0)\ne 0,$ and it is a fold if, in addition,
one has $P(\xi)\ne (0,0)^\top,$ where $\xi$ is the tangent vector to the singular curve through $(0,0)$. When $g_{xx}(0,0)\ne 0$,
the singular curve has a representation $(l(v),v),$ $l(0)=0,$ so $\xi$ is $(l'(0),1)^\top.$ Since $l$ solves the equation $g_x(x,v)=0,$ then
$l'(0)=-g_{xv}(0,0)/g_{xx}(0,0)$ and $P(\xi)=(g_v(0,0),1)^\top\ne (0,0)^\top.$
Thus, the singular point $(0,0)$ is indeed the fold, if $g_x(0,0)=0$ but $g_{xx}(0,0)\ne 0$.

Now assume $g_x(0,0)=g_{xx}(0,0)=0$ but $g_{xv}(0,0)\ne 0.$ Then the singular curve has a representation
$(x,r(x)),$ $r(0)=0,$ and vector $\xi$ is $(1,r'(0))^\top$. Because $r(x)$ again solves the equation $g_x(x,v)=0,$ we have the equality
$$
r'(0)=-\frac{g_{xx}(0,0)}{g_{xv}(0,0)}=0.
$$
It follows from here that $P(\xi)= (0,0)^\top.$ Therefore, if $g_{xx}(0,0)=0$ the singular point $(0,0)$ is not a fold.
In this latter case to verify it to be a cusp one needs to calculate $d^2(r(x),g(x,r(x)))/dx^2$ at the point $(0,0).$ The calculation
gives the formula
$$
(g_{xx}(0,0)+2g_{xv}(0,0)r'(0)+g_{vv}(r'(0))^2 + g_v(0,0)r^{\prime\prime}(0),r^{\prime\prime}(0)) =
(g_v(0,0)r^{\prime\prime}(0),r^{\prime\prime}(0)).
$$
due to equalities $g_{xx}(0,0)=0,$ $r'(0)=0.$ Thus, if $r^{\prime\prime}(0)\ne 0$ the second derivative is not zero vector, and the point is a cusp.
For the derivative $r^{\prime\prime}(0)$ we have
\begin{equation}
r^{\prime\prime}(0)= -\frac{g_{xxx}(0,0)}{g_{xv}(0,0)}.
\end{equation}
Thus, the conditions for a singular point to be a cusp casts in two equalities and two
inequalities
\begin{equation}\label{cusp}
g_x(0,0)= g_{xx}(0,0)=0,\;g_{xv}(0,0)\ne 0 \;g_{xxx}(0,0)\ne 0.
\end{equation}

A position of the singular curve w.r.t. levels $H=c$ is important. In particular,
we need to know when the intersection of this curve and submanifold $H=H(s)$
is transverse and when is not the case. For the transverse case for all $c$ close enough to $H(s)$ the intersection of
this curve with levels $H=c$ will be also transversal and hence, consists of one
point. For the nontransverse case we need to know what will happen for a
close levels. In a sense, the transversality condition can be considered
as some genericity condition for the chosen function $H$.

The point $s$ is a fold point for the restricted map $p_r : SM \to B$, if, in addition to
the equality $\Delta(s) = 0,$ i.e. $g_{x}(0,0)=0$, the condition of quadratic
tangency $g_{xx}(0,0)\ne 0$ will be satisfied. This is equivalent to the inequality $\Delta_{x}(s) \ne 0$ and allows
one to express from the equation $g_x =0$ the singular curve on $SM$  near $s$ in the form $x=l(v),$
$l(0)=0.$ To determine if this curve intersects transversely the level
$H=H(s)$, let us calculate the derivative
$$
\frac{d}{dv}H(l(v),f(l(v),v),g(l(v),v),v)|_{v=0}= H_u(0,0,0,0)g_v(0,0)+H_v(0,0,0,0).
$$
Taking into account that $f(x,v),g(x,v)$ are solutions of the system $H_y =
0,$ $H_x = 0$ near point $s$ we сan calculate $g_v(0,0)$, then the numerator of this derivative at point
$s$ does not vanish if
\begin{equation}\label{parab}
H_{xy}[H_uH_{yv}-H_vH_{yu}]-H_{yy}[H_uH_{xv}-H_vH_{xu}] \ne 0.
\end{equation}
It is hard to calculate this quantity when $H$ is taken in a general
form. Therefore we transform $H$ near $s$ to a more tractable form w.r.t. variables
$(x,y)$ with parameters $(u,v).$ In order not to care about the smoothness we assume henceforth all functions are $C^\infty$.
The following assertion is valid
\begin{lemma}\label{y_square}
Suppose a smooth function $H(x,y,u,v)$ be given such that at parameters
$(u,v)=(0,0)$ the Hamiltonian system
$$\dot x = H_y,\;\dot y = - H_x$$
has a degenerate equilibrium $(x,y)=(0,0)$ of the linearized system
with the double zero non-semisimple eigenvalue (the related Jordan form is two dimensional).
Then there exists a $C^{\infty}$-smooth transformation $\Phi: (x,y)\to (X,Y)$ smoothly depending on parameters $(u,v)$ and
respecting 2-form $dx\wedge dy$ such that the Hamiltonian $H\circ \Phi$ in new variables $(X,Y)$ takes the form
\begin{equation}\label{fold}
H(X,Y,u,v)= h(X,u,v)+ H_{1}(X,Y,u,v)Y^2,
\end{equation}
\end{lemma}
\proof We act as follows.
Consider first the Hessian $(a_{ij})$ of $H$ in variables $(x,y)$ at the
point $(x,y)=(0,0)$ on the leaf $(u,v)=(0,0).$ Its determinant vanishes
but not all its entries are zeroes and its rank is 1, since we supposed the eigenvalues to be non-semisimple ones.
Then one has either $a_{11}\ne 0$ or
$a_{22}\ne 0$ due to symmetricity of the Hessian and the assumption on non-semisimplicity of zero eigenvalue. We assume $a_{22}\ne 0$
to be definite, for our case above it follows from the assumption
that minor $H_{yy}H_{xu}-H_{xy}H_{yu}\ne 0.$ We first solve the equation $H_y =0$ in a neighborhood of point $(0,0,0,0).$
In view of the implicit function theorem, since $a_{22}= H_{yy}(0,0,0,0)\ne 0$, the equation has a solution $y=f(x,u,v),$
$f(0,0,0)=0.$ After shift transformation $x=X,$
$y= Y + f(x,u,v)$ we get the transformed Hamiltonian $\hat H(X,Y,u,v)$ which is of the form
$$
\begin{array}{l}
\hat H(X,Y,u,v)=h(X,u,v)+ Y \tilde{H}(X,Y,u,v),\;\tilde{H}(X,0,u,v)= H_y(x,f(x,u,v),u,v)\equiv 0,\\h(X,u,v)=H(x,f(x,u,v),u,v).
\end{array}
$$
Thus function $\tilde{H}$ can be also represented as
$\tilde{H}= Y\bar{H},$ $\bar{H}(0,0,0,0)= H_{yy}(0,0,0,0)/2 \ne 0.$ It
is worth noticed that $(u,v)$ are considered as parameters, therefore the
transformation $(x,y)\to (X,Y)$ respects the 2-form: $dx\wedge dy = dX\wedge dY$.$\blacksquare$

We now restore the notations $(x,y)$ and assume $H$ being in the form
(\ref{fold}). Then the condition $(0,0,0,0)$ to be the equilibrium of
the fast system is $y=0,$ $h_x(0,0,0)=0.$ The claim it to be
degenerate and non-semisimple is $h_{xx}(0,0,0)=0,$ since $H_1(0,0,0,0)\ne
0$. In this case the requirement the slow manifold to be a smooth solution of the equation
$h_x(x,u,v)=0$ near $(0,0,0,0)$ leads to inequalities $h_{xu}(0,0,0)\ne 0$
or $h_{xv}(0,0,0)\ne 0$. This can be always converted into the first one renaming slow variables, we assume this to be the case.
The latter means that the slow manifold has a representation $y=0,$
$u=g(x,v),$ $g(0,0)=0,$ $g_x(0,0)=0.$ At last, the point $(0,0,0,0)$ on
the slow manifold will be a fold, if $g_{xx}(0,0)\ne 0$, that is equivalent to inequality $h_{xxx}(0,0,0)\ne 0.$

Now we expand $h(x,u,v)$ in $x$ up to the third order terms
$$
H(x,y,u,v)=h(x,u,v)+  H_1(x,y,u,v)y^2 = h_0(u,v)+a(u,v)x+b(u,v)x^2 + c(u,v)x^3 + O(x^4) + H_1y^2.
$$
Here one has $a(0,0)=0$, $a_u(0,0)\ne 0,$ $b(0,0)=0,$ $c(0,0)\ne 0.$ We
have some freedom to change parameters $(u,v).$ Using inequality $a_u(0,0)\ne
0,$ we introduce new parameters $u_1 = a(u,v).$ To get $v_1$ we apply a
symplectic transformation using 2-form $du\wedge dv$. To that end we
express $u=\hat{a}(u_1,v)= R_{v}$, where for $|u_1|, |v|$ small enough
$$
R(u_1,v)=\int\limits_0^v \hat{a}(u_1,z)dz,\;\frac{\partial^2 R}{\partial u_1\partial v} = \frac{\partial \hat a}{\partial u_1}\ne 0.
$$
Then one has $v_1 = R_{u_1}$ and $du\wedge dv = du_1\wedge dv_1.$ After
this transformation which does not touch variables $x,y$ we come to the
form of $H$
\begin{equation}\label{param}
H(x,y,u_1,v_1)= h_0(u_1,v_1)+ u_1x+\hat{b}(u_1,v_1)x^2 + \hat{c}(u_1,v_1)x^3 + O(x^4) + \hat{H}_1 y^2.
\end{equation}
In this form we can check the transversality of the singular curve on $SM$ and submanifold $H=h_0(0,0)=c_0$ at the point $(0,0,0,0)$.
Since we have $H^0_{xy}=0$, $H^0_{yy}\ne 0$ (zeroth upper script means the calculation at the point
$(0,0,0,0)$), then the inequality (\ref{parab}) casts (we restore notations $u,v$ again) as
$$
H^0_u H^0_{xv} - H^0_v H^0_{xu}\ne 0,
$$
that is expressed via $h_0$ and $a = u$ looks as follows
$$
-\frac{\partial h_0}{\partial v}(0,0)\ne 0.
$$
Thus, we come to the conclusion:

{\em if a function $H$ is generic and the singular point $s$ on $SM$ is a fold for the mapping $p_r$, this is
equivalent to the condition that this point on the related symplectic leaf is
parabolic and the unfolding of $H$ in parameters $(u,v)$ is generic.}

At the next step we want to reduce the dimension of the system near the
fold point $s\in SM$ and get a smooth family of nonautonomous Hamiltonian systems in one degree
of freedom. This will allow us to describe the principal part of the system near
singularity using some rescaling for the system near $s$. We
shall work in coordinates where $H$ takes the form (\ref{fold}).

In Darboux-Weinstein coordinates $(x,y,u,v)$ the slow-fast Hamiltonian system with Hamiltonian $H$ (\ref{fold}) is written near $s$ as follows
\begin{equation}\label{sfloc}
\dot x = H_y,\;\dot y = - H_x,\;\dot u = \eps H_v,\;\dot v = - \eps H_u.
\end{equation}
Without a loss of generality, one can assume $H(s)=0$. Since $H_v(s)\ne 0$, then near $s$ levels $H=c$ for $c$ close to zero are given as graphs of
a function $v=S(x,y,u,c),$ $S(0,0,0,0)=0$, $S_c = 1/H_v \ne 0$. These graphs intersect
transversely the singular curve near $s$, thus it occurs at one point on the related graph. The
intersection of $SM$ with a level $H=c$ is a smooth curve (the slow curve in this level) with unique
tangency point with the leaves $(u,v)=(u_0,v_0)$ at the related point of the singular
curve.

Let us perform the isoenergetic reduction of system (\ref{sfloc}) on the level $H=c$, then $S$ be a new (nonautonomous)
Hamiltonian and $u$ be the new ``time''. After reduction the system transforms
\begin{equation}\label{sfnon}
\eps\frac{dx}{du} = H_y/H_v = - S_y,\;\eps\frac{dy}{du} = - H_x/H_v = S_x.
\end{equation}
If we introduce fast time $\tau = u/\eps,$ we come to the nonautonomous Hamiltonian system which depends on slowly varying time $\eps\tau$.

A slow curve for this system is given by equations $S_y =0,$ $S_x =0$ that is led to the same system $H_y =0,$ $H_x =0$
where one needs to plug $S$ into $H$ instead of $v$. Thus we obtain the intersection of $SM$ with the level
$H=c$. For the case of $H$ we consider, this curve has a representation $y=0,$ $u-u_c = a(c)(x-x_c)^2 +
o((x-x_c)^2),$ $a(0)\ne 0,$ with $(x_c,0,u_c)$ being the coordinates of the singular curve trace on the level $H=c$ for the fixed $c$ close to
$c=0$.  This is extracted from the system $u=g(x,v)$, $v=S(x,0,u,c)$, therefore $u$ is
expressed via $x$ by the implicit function theorem, since $1-g_vS_u = (H_u g_v + H_v)/H_v \ne 0$ at the point
$(x,y,u,c)=$ $(0,0,0,0)$ (and any point close to it on the singular
curve).

Fast or frozen system for the reduced system is given
by setting $u=u_0$ and varying $u_0$ near $u_0 = u_c.$ On the leaf $u_0=u_c$ we get a one degree of freedom
Hamiltonian system with the equilibrium at $(x_c,0)$ having double zero non-semisimple eigenvalue.
Using the form (\ref{fold}) of Hamiltonian we come to the
one-degree-of-freedom system
\begin{equation}\label{nonaut_fold}
\frac{dx}{d\tau}=\frac{2yH_1 + y^2H_{1y}}{h_{0v}+b_v x^2 + c_v x^3 +y^2H_{1v}+ O(x^4)},\;
\frac{dy}{d\tau}= -\frac{\eps\tau + 2b x + 3c x^2 +y^2H_{1x}+O(x^3)}{h_{0v}+b_v x^2 + c_v x^3 +y^2H_{1v}+ O(x^4)}.
\end{equation}

We can also find the foliation of $SM$ near point $s$ into level lines of
Hamiltonian $H$ restricted on this manifold. It is nothing else as the local phase portrait of a slow system near a singular curve.
It is given by levels of the function
$\hat h = H(x,0, g(x,v),v)=$ $h_0(g(x,v),v)+  g(x,v)x + b(g(x,v),v)x +$ $c(g(x,v),v)x^3$.
The manifold $SM$ has the line of folds (= the singular curve).
This line is projected in $B$ near $s$ as a smooth curve (this curve can be called as a {\em discriminant curve}
similar to the theory of implicit differential equations, see \cite{DopGl}) in such a way that the image of $SM$
lies from the one side of this curve. Then foliation curves are projected as having cusps at the discriminant
curve. This picture is the same as in a nonhamiltonian slow-fast system
with one fast and two slow variables (see \cite{Arn}).

\subsection{Rescaling near a fold}

Here we want to find a principal part of the system (\ref{nonaut_fold})
near a fold point. We use a blow up method like in \cite{DR,Kr_Sz}. It is
not surprising that we meet the Painlev\'e-1 equation here, see
\cite{Haberman}. We consider our derivation as a more direct and
consequent one.

We start with the observation that the one degree of freedom nonautonomous
Hamiltonian $S(x,y,u,c)$ can be written in the form
(\ref{param}), if we expand it near the point $(x_c,0,u_c)$ being the
trace of the singular curve for $H$
$$
S(x,y,u,c)= s_0(u,c)+\alpha(u,c)(x-x_c)+\beta(u,c)(x-x_c)^2+\gamma(u,c)(x-x_c)^3 + O((x-x_c)^4)+y^2S_1(x,y,u,c),
$$
where $\alpha(u_c,c)=\beta(u_c,c)=0,$ $\alpha_u(u_c,c)\ne 0,$ $\gamma(u_c,c)\ne
0,$ $S_1(x_c,0,u_c,c)\ne 0.$ Denote $x-x_c = \xi.$ This gives us the following form of the reduced system
\begin{equation}
\begin{array}{l}
\displaystyle{\eps\frac{d\xi}{du}= - \frac{\partial S}{\partial y}= -2yS_1 - y^2\frac{\partial S_1}{\partial y}},\\
\\\displaystyle{\eps\frac{dy}{du}= \frac{\partial S}{\partial \xi}= \alpha(u,c)+ 2\beta(u,c)\xi + 3\gamma(u,c)\xi^2 + O(\xi^3)+
y^2\frac{\partial S_1}{\partial \xi}}.
\end{array}
\end{equation}
Now we use the variable $\tau = (u-u_c)/\eps$ and add two more equations $u' = du/d\tau = \eps$ and $\eps' = 0$ to the system.
Then the suspended autonomous system will have an equilibrium
at the point $(\xi,y,u,\eps)=(0,0,u_c,0)$.
The linearization of the system at this equilibrium has a matrix being nothing
else as 4-dimensional Jordan box. To study the solutions to this system near the equilibrium we,
following the idea in \cite{DR,Kr_Sz} (see also a close situation in \cite{chiba}), blow up a neighborhood of this point
by means of the coordinate change
\begin{equation}\label{blow}
\xi=r^2X,\;y=r^3Y,\;u=r^4Z,\;\eps=r^5,\; r\ge 0.
\end{equation}
Since $\dot \eps = 0$ we consider $r = \eps^{1/5}$ as a small
parameter. The system in these variables takes the form
$$
X' = -r(2YS_1(x_c,0,u_c,c)+\cdots],\;\dot Y = r[\alpha_u(u_c,c)Z +3\gamma(u_c,c)X^2 + \cdots],\;\dot Z = r.
$$
After re-scaling time $r\tau = s$, denoting $' = d/ds$, setting $r = 0$ we get
$$
X' = -2Y,\;Y' = \alpha_c Z +3\gamma_c X^2,\;Z' = 1.
$$
where $\alpha_c = \alpha_u(u_c,c),$ $\gamma_c = \gamma(u_c,c).$ The obtained system
 is equivalent to the well known Painlev\'e-I equation \cite{Pain,Haberman,Novok}
$$\frac{d^2 X}{dZ^2} + 2\alpha_0Z + 6\gamma_0 X^2
=0.$$
By scaling variables this equation can be transformed to the standard
form.
It is known \cite{Haberman} that this equation appears
at the passage through a parabolic equilibrium in the fast system. We come
to this equation directly using blow-up procedure.

In fact, it is not the all story if one wants to study in details the behavior near the fold
point: one need to derive more six systems. In fact the blow-up procedure
means that we blow up the singular point of the suspended 4-dimensional vector field
with $\dot \eps =0$ added till the 3-dimensional sphere $S^3$ and a
neighborhood of the singular point do till a neighborhood of this sphere
with $r\ge 0$ being a coordinate in the transverse direction to the
sphere. In order to get a sphere one needs to regard $X^2+Y^2+Z^2+E^2 =1$.
But it is not convenient to work in coordinates on the sphere, therefore
it is better to work in charts. These charts are obtained, if we
consequently set $E=1$, then we get the system derived above. Here we consider only the chart with $E=1$, but not $E=-1$, because
we assume $\eps >0$. Other six charts are obtained if we set $X=\pm 1,$ $Y+\pm
1$ and $Z=\pm 1$, then other six system will be derived. In the principal
approximation one needs to take limits $r\to 0$ in this systems. The
recalculation from one system to another one is given by the main formulae
for the blow-up. Then one can understand the whole picture of a passage of
orbits through the neighborhood of the disruption point. This will be done
elsewhere.

As is known, for 2-dim slow-fast (dissipative) system such the passage is described by the
Riccati equation \cite{MR,Kr_Sz}. Here we get the Painlev\'e-I equation.

\section{Cusp for the slow manifold projection}

For the case when $s$ is a cusp, the related singular curve on $SM$ is tangent to the symplectic leaf through $s$
(see above the equality $P(\xi)=0$). Due to the last inequality in (\ref{cusp}), this tangency takes place at only point
$s$, other points of the singular curve near $s$ are folds. Below without a loss of generality we assume $s$ to be the origin $(0,0,0,0).$

The singular curve is also tangent at $s$ to the level $H=H(s)$. Indeed, for the case of a cusp the singular curve on
$SM$ has a representation $x=r(v),$ $r(0)=r'(0)=0$ near $s$ (see above). Then from equalities $H_x(0,0,0,0) = H_y(0,0,0,0) = 0,$
$g_x(0,0),$ $r'(0)= 0$ the tangency follows
$$
\frac{d}{dx}H(x,f(x,r(x)),g(x,r(x)),r(x))|_{x=0}=0.
$$

Consider for $\eps =0$ the leaf $F_b,$ $b=p(s),$ of the symplectic foliation and the
Hamiltonian system with Hamiltonian $H$ restricted to this leaf (the fast Hamiltonian system). This
one degree of freedom system has an equilibrium at $s$. This equilibrium is degenerate: it has double zero eigenvalue
as for a fold but this equilibrium is more degenerate than that being parabolic one. To make further calculations we again use
lemma \ref{y_square} and the normal form method in some neighborhood of this equilibrium. We want to show the equilibrium
be of co-dimension 2 though it has the same linear part as for a parabolic point, but it obeys two additional equalities. The partial normal
form for such an equilibrium depending on parameters $(u,v)$ looks as follows
\begin{equation}\label{cusp_ham}
H(x,y,u,v)= h_0(u,v)+ a_1(u,v)x + \frac{a_2(u,v)}{2}x^2 + \frac{a_3(u,v)}{3}x^3 + \frac{a_4(u,v)}{4}x^4 + O(x^5)+ y^2H_1(x,y,u,v)
\end{equation}
with $dh_0 (0,0)\ne 0, $ $H_1(0,0,0,0)\ne 0$, $a_4(0,0)\ne 0,$ $a_1(0,0)=0$. The first condition for the co-dimension 2 here
is equalities
$a_2(0,0)=a_3(0,0)=0.$ These two conditions follow from the assumption for $s$ to be a cusp for the map $p: M\to B$,
then $g_x(0,0)=0$ implies $a_2(0,0)=0$ and $g_{xx}(0,0)=0$ does $a_3(0,0)=0.$ Inequality $g_{xxx}(0,0)\ne 0$ mean here $a_4(0,0)\ne 0$.
It appears the sign of $a_4(0,0)$ is essential and opposite
signs lead to different structures of nearby fast systems on the
neighboring leaves (see Figs. 1a-1b).
\begin{figure}[h]
\centering\includegraphics[width=0.6\textwidth]{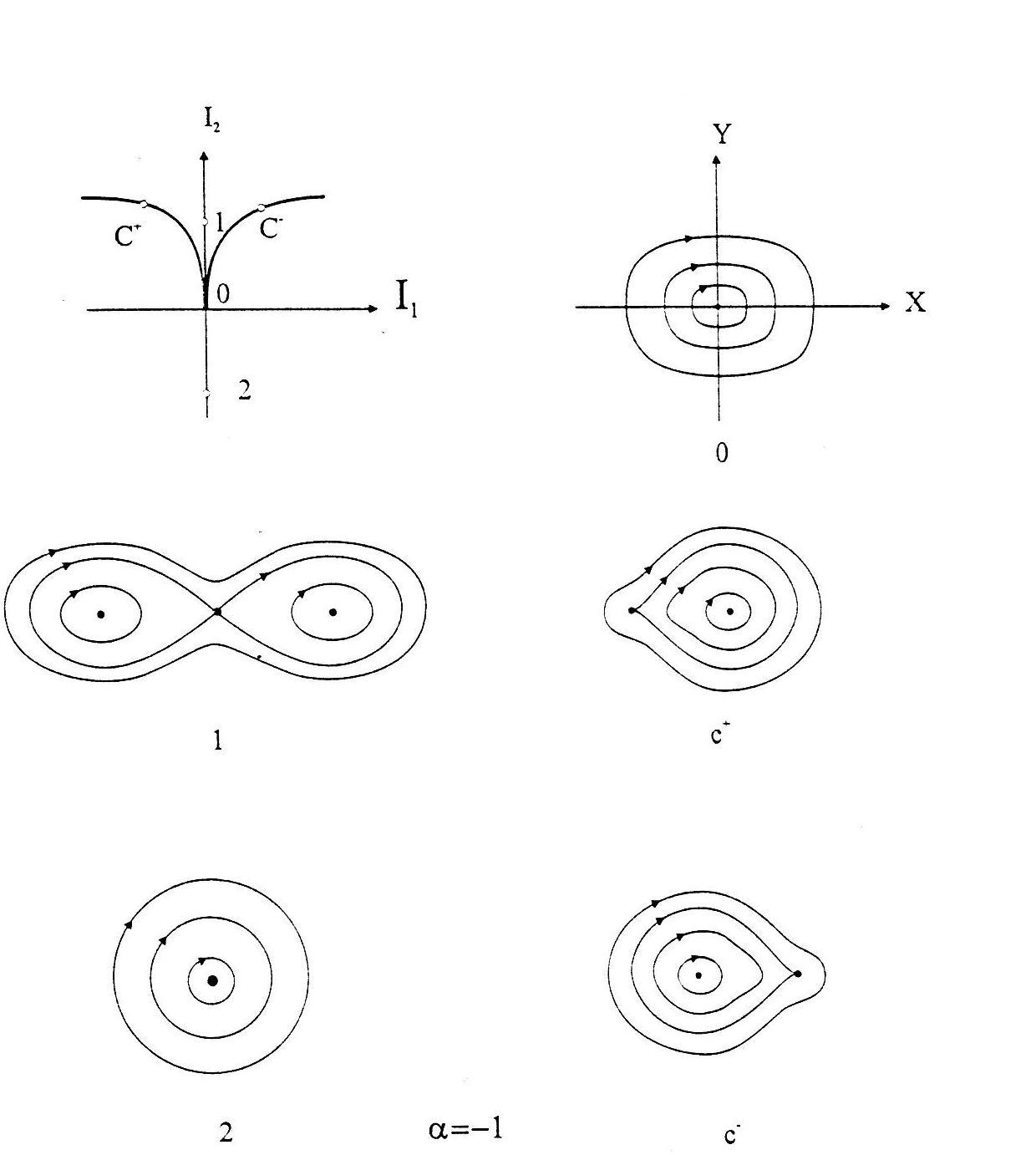}
\end{figure}

\begin{figure}[h]
\centering\includegraphics[width=0.6\textwidth]{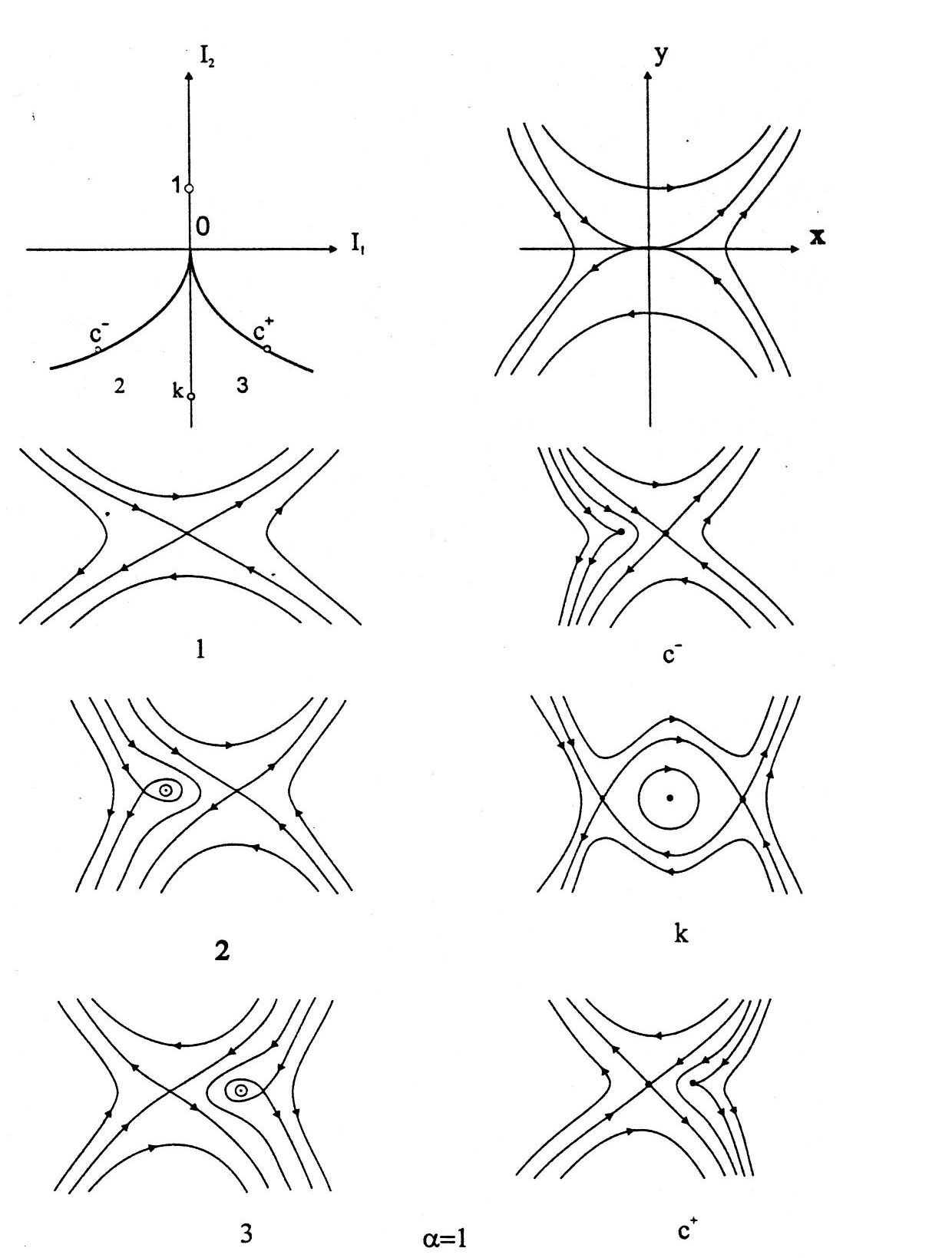}
\end{figure}

In order the unfolding in $(u,v)$ would be generic, the
condition ${\rm det} (D(a_1,a_2 )/D(u,v))\ne 0$ at $(u,v)=(0,0)$ has to be met. This
just corresponds to that when rank of $(3\times 4)$-matrix above is 3 and
then the curve of singular points here is expressed as
$(x,y(x),u(x),v(x))$ and $y'(0)=u'(0)=$ $v'(0)=0.$ The projection of the
singular curve on $B$ ($(u,v)$-coordinates) is the cusp-shaped curve which
is up to higher order terms is given parametrically in $x$ as
$$a_1(u,v)+a_2(u,v)x + a_3(u,v)x^2 + a_4(u,v)x^3 +O(x^4)=0,\; a_2(u,v) + 2a_3(u,v)x + 3a_4(u,v)x^2 +O(x^3)=0.$$

To ease again the further calculations we take functions $a_1,a_2$ as new parameters instead
of $u,v$ using the fact that ${\rm det} (D(a_1,a_2)/D(u,v))\ne 0$.
Keeping in mind that $(u,v)$ are symplectic coordinates w.r.t. 2-form $du\wedge
dv$ we make the change via a symplectic transformation. To this end, we
assume first the determinant be positive, otherwise we make the
redesignation $(a_1,a_2)\to (-a_1,a_2).$ One of partial derivatives of $a_1$ in
$u$, $v$ at $(0,0)$ does not vanish, so one can take $a_1$ as a new slow variable $u_1 = a_1(u,v).$
Adding to that some $v_1$ in order to get a symplectic transformation $(u,v)\to (u_1,v_1)$ we
come to the same form of $H$ w.r.t the new variables $(u_1,v_1)$ and new coefficients $a_1(u_1,v_1)=u_1, a_2(u_1,v_1)$ which we
again denote as $(u,v)$ and $u,a_2(u,v)$. Then we get $\partial a_2/\partial v \ne 0$, since ${\rm det} D(a_1,a_2)/D(u,v)\ne 0$.

At the next step we want to use the isoenergetic reduction and get again a
family of nonautonomous Hamiltonian system in one degrees of freedom depending on a parameter $c$ in a neighborhood of point $s=(0,0,0,0)$.
To this end, we have to know which of partial derivatives of $h_0$ is
nonzero. Recall that the singular curve on $SM$ at the point $s$ is
tangent  to $H=H(s)$. We assume without loss of generality that $H(s)=0$.
For Hamiltonian (\ref{cusp_ham}) the derivative of $H$ along the tangent
vector to singular point at $s$ is equal to $H^0_u H^0_{xv}-H^0_v H^0_{xu}
=0$ but $H^0_{xv}= 0$ and $H^0_{xu} = 1$ (since $a_1(u,v)=u$), hence one
has $H^0_v =0.$ Due to the assumption $dH \ne 0$ at $s$ we come to the inequality $H^0_{u} \ne 0.$
This implies the equation $H=c$ near $s$ be solved as $u=S(x,y,v,c),$
$S(0,0,0,0)=0.$ The derivative $S^0_c$ does not vanish, since it is equal
at $s$ to $(\partial h/\partial u)^{-1} \ne 0.$
\begin{equation}\label{nonaut_cusp}
\begin{array}{l}
\displaystyle{\frac{dx}{d\tau}=-\frac{2yH_1 + y^2H_{1y}}{h_{0u}+ x + a_{2u}x^2/2 + a_{3u} x^3/3 + a_{4u} x^4/4 + O(x^5)+ y^2H_{1u}}},\\
\displaystyle{\frac{dy}{d\tau}= \frac{u + a_2 x + a_3 x^2 + a_4x^3 + y^2H_{1x}+O(x^4)}{h_{0u}+ x + a_{2u}x^2/2 + a_{3u} x^3/3 +
a_{4u} x^4/4 + O(x^5)+ y^2H_{1u}}},\\
\displaystyle{\frac{dv}{d\tau}=\eps,\;\frac{d\eps}{d\tau}}=0.
\end{array}
\end{equation}
We added here more two equations $dv/d\tau =\eps$ and $d\eps/d\tau =0$ and get again an equilibrium at
$(0,0,0,0).$ To study the system near the equilibrium we perform the blow-up transformation
\[
x\to rX,\;y\to r^2Y,\;v\to r^2Z,\;u\to r^3C\;\eps \to r^3E.
\]
We present here only the system with $E=1$, then $r=\eps^{1/3}.$ After writing the system in new coordinates,
scaling time $r\tau =s$, setting $r=0$ and denoting constants $\rho = 2H_1(0,0,0,0)(\partial h_0(0,0)/\partial u)^{-1},$
$\sigma = (\partial h_0(0,0)/\partial u)^{-1},$ $\beta = (\partial a_{2u}(0,0)/\partial v)(\partial h_0(0,0)/\partial u)^{-1}$,
$\alpha=a_4(0,0),$ $A=C(\partial h_0(0,0)/\partial u)^{-1}$, we come to the system
\[
\dot X = - \sigma Y,\;\dot Y = A + \beta ZX + \alpha X^3,\;\dot Z = 1,
\]
that is just the Painlev\'e-II equation. If one introduce the new time $s\sigma = \xi$ and reduce it to the standard form.

It represents the behavior of solutions inside of a 4-dimension disk.
Other equations should be written and understood in order to catch the whole picture.

Thus, we have proved a theorem that gives a connection between an orbit
behavior of a slow fast Hamiltonian system near its disruption point and
solutions to the related Painlev\'e equations. We formulate these results
as follows.
\begin{theorem}
Let a smooth slow fast Hamiltonian vector field with a Hamiltonian $H$ be given on a smooth
Poisson bundle $p: M\to B$ with a Poisson 2-form $\omega$ and $B$ being a smooth symplectic manifold with a
symplectic 2-form $\nu$. We endow $M$ with the singular symplectic structure $\omega
+\eps^{-1}p_*\nu$. Suppose the set of zeroes $SM$ of the fast vector
fields on symplectic leaves generated by $H$ forms a smooth submanifold in
$M$. Points of $SM$ being tangent to symplectic leaves and other than critical points of $H$  we call to be
disruption points. If $M$ is four-dimensional and $B$ is two-dimensional
ones then generically a disruption point $s$ can be only of two types w.r.t. the map $p|_{SM}$: a
fold or a cusp. Then the slow-fast system near $s$ can be reduced in the principal approximation
after the isoenergetical reduction and some blow-up to either
the Painlev\'e-I equation, if $s$ is a fold, or to the Painlev\'e-II
equation, if $s$ is a cusp.
\end{theorem}

\section{Acknowledgement}

L.M.L. acknowledges a partial support from the Russian Foundation for Basic Research (grant 14-01-00344a)
and the Russian Science Foundation (project 14-41-00044).
E.I.Y. is thankful for a support to the Russian Ministry of Science and Education
(project 1.1410.2014/K, target part).

\end{document}